\title{On two generalizations of the Alon-Tarsi polynomial method}
\author{{Dan Hefetz \thanks{Institute of Theoretical Computer Science, ETH Zurich,
CH-8092 Switzerland. Email: dan.hefetz@inf.ethz.ch.}}}

\documentclass[12pt]{article}
\usepackage{amsmath,amssymb,latexsym,color,epsfig,a4}

\newif\ifnotesw\noteswtrue

\def\({\left(}
\def\){\right)}

\newtheorem{theorem}{Theorem}[section]
\newtheorem{lemma}[theorem]{Lemma}

\newtheorem{definition}[theorem]{Definition}

\newtheorem{observation}[theorem]{Observation}
\newtheorem{corollary}[theorem]{Corollary}
\newtheorem{conjecture}[theorem]{Conjecture}

\renewcommand{\epsilon}{\varepsilon}

\begin{document}
\maketitle

\begin{abstract}
In a seminal paper~\cite{AT92}, Alon and Tarsi have introduced an algebraic technique for proving upper bounds on the choice number of graphs (and thus, in particular, upper bounds on their chromatic number). The upper bound on the choice number of $G$ obtained via their method, was later coined the \emph{Alon-Tarsi number of $G$} and was denoted by $AT(G)$ (see e.g.~\cite{Toft}). They have provided a combinatorial interpretation of this parameter in terms of the eulerian subdigraphs of an appropriate orientation of $G$. Their characterization can be restated as follows. Let $D$ be an orientation of $G$. Assign a weight $\omega_D(H)$ to every subdigraph $H$ of $D$: if $H \subseteq D$ is eulerian, then $\omega_D(H) = (-1)^{e(H)}$, otherwise $\omega_D(H) = 0$. Alon and Tarsi proved that $AT(G) \leq k$ if and only if there exists an orientation $D$ of $G$ in which the out-degree of every vertex is strictly less than $k$, and moreover $\sum_{H \subseteq D} \omega_D(H) \neq 0$. Shortly afterwards~\cite{Alon}, for the special case of line graphs of $d$-regular $d$-edge-colorable graphs, Alon gave another interpretation of $AT(G)$, this time in terms of the signed $d$-colorings of the line graph. In this paper we generalize both results. The first characterization is generalized by showing that there is an infinite family of weight functions (which includes the one considered by Alon and Tarsi), each of which can be used to characterize $AT(G)$. The second characterization is generalized to all graphs (in fact the result is even more general --- in particular it applies to hypergraphs). We then use the second generalization to prove that $\chi(G) = ch(G) = AT(G)$ holds for certain families of graphs $G$. Some of these results generalize certain known choosability results.
\end{abstract}

\section{Introduction}\label{sec::intro}
Let $G=(V,E)$ be a graph and let $k$ be a positive integer. A \emph{proper $k$-coloring} of $G$ is a mapping $c : V \rightarrow \{0,1, \ldots, k-1\}$ such that $c(u) \neq c(v)$ whenever $(u,v) \in E$. The smallest integer $k$ for which $G$ admits a proper $k$-coloring is called the \emph{chromatic number} of $G$ and is denoted by $\chi(G)$. A \emph{proper list coloring} of $G$ from lists $\{L(u) : u \in V\}$ is a mapping $f : V \rightarrow \bigcup_{u \in V} L(u)$ such that $f(u) \neq f(v)$ whenever $(u,v) \in E$ and $f(u) \in L(u)$ for every $u \in V$. The smallest integer $k$ for which $G$ admits a proper list coloring from any set of lists $\{L(u) : u \in V\}$ such that $|L(u)| \geq k$ for every $u \in V$, is called the \emph{choice number} of $G$ and is denoted by $ch(G)$. The choice number was introduced by Vizing~\cite{Vizing} and independently by Erd\H{o}s et al.~\cite{ERT}.

Clearly, by definition, $ch(G) \geq \chi(G)$ for every graph $G$; however the converse inequality does not hold in general. Indeed, it was proved in~\cite{ERT} that $ch(K_{n,n}) = (1+o(1))\log_2 n$ and so, in general, the choice number of a graph cannot be bounded from above by any function of its chromatic number. A graph which satisfies $ch(G) = \chi(G)$ is called \emph{chromatic-choosable}. A lot of effort has been invested in proving that certain graph families are chromatic-choosable, and, in particular, many conjectures hypothesizing sufficient conditions for a graph to be chromatic-choosable were made. Probably, the most famous of these is the following (see e.g.~\cite{Toft}):

\begin{conjecture} [List Edge Coloring Conjecture] \label{LCC}
Every line-graph is chromatic-choosable.
\end{conjecture}

Despite many efforts by various researchers (e.g.~\cite{Galvin, Kahn} to name just a few), Conjecture~\ref{LCC} is still open. 

\subsection{The Alon-Tarsi polynomial method} \label{subsec::ATmethod}

Let $G=(V,E)$ be an undirected multi-graph with vertex set $\{1, \ldots, n\}$. The \emph{graph polynomial of $G$} is defined by 
$$
P_G(x_1, x_2, \ldots, x_n) = \prod_{1 \leq i < j \leq n \atop (i,j) \in E}(x_i - x_j).
$$
If $E = \emptyset$, then we define $P_G \equiv 1$.

This polynomial was introduced by Petersen~\cite{Petersen} more than a century ago. It has been extensively studied since by various researchers (see e.g.~\cite{AT97, LiLi, Lovasz, Tarsi}).

It is clear that the graph polynomial encodes information about its proper colorings. Indeed, a graph $G$ is $k$-colorable if and only if there exists an $n$-tuple $(a_1, a_2, \ldots, a_n) \in \{0,1, \ldots, k-1\}^n$ such that $P_G(a_1, a_2, \ldots, a_n) \neq 0$. Similarly, $G$ is $k$-choosable if and only if for every family of sets $\{S_i \subseteq \mathbb{R}: 1 \leq i \leq n\}$, each of size at least $k$, there exists an $n$-tuple $(a_1, a_2, \ldots, a_n) \in S_1 \times S_2 \times \ldots \times S_n$ such that $P_G(a_1, a_2, \ldots, a_n) \neq 0$. 

The following theorem gives a sufficient condition for the existence of such an $n$-tuple.

\begin{theorem}(\textbf{Combinatorial Nullstellensatz}~\cite{CNSS}) \label{t::CNSS}
\label{Combnull} Let $\mathbb{F}$ be an arbitrary field, and let
$f=f(x_{1},\ldots,x_{n})$ be a polynomial in
$\mathbb{F}[x_{1},\ldots,x_{n}]$. Suppose that the degree $deg(f)$ of
$f$ is $\sum_{i=1}^{n}t_{i}$, where each $t_{i}$ is a non-negative
integer, and suppose that the coefficient of
$\prod_{i=1}^{n}x_{i}^{t_{i}}$ in $f$ is nonzero. Then, if
$S_{1},\ldots,S_{n}$ are subsets of $\mathbb{F}$ with $|S_{i}| >
t_{i}$, then there are $s_{1} \in S_{1}, s_{2} \in S_{2},\ldots,s_{n}
\in S_{n}$ so that $f(s_{1},\ldots,s_{n}) \neq 0$.
\end{theorem}

We are now ready to define the main concept of this paper.

\begin{definition} \label{algChoosable}
Let $G=(V,E)$ be a multi-graph with vertex set $V = \{v_1, \ldots, v_n\}$ and let $f : V \rightarrow \mathbb{N}^+$. We say that $G$ is \emph{Alon-Tarsi $f$-choosable}, or $f-AT$ for brevity, if there exists a monomial $c \prod_{i=1}^n x_i^{t_i}$ in the expansion of $P_G$ such that $c \neq 0$ (say, in $\mathbb{R}$) and $t_i \leq f(i)-1$ for every $1 \leq i \leq n$. If $f(i)=k$ for every $1 \leq i \leq n$, then we say that $G$ is \emph{Alon-Tarsi $k$-choosable}, or $k-AT$ for brevity. The smallest integer $k$ for which $G$ is Alon-Tarsi $k$-choosable, denoted by $AT(G)$, is called the \emph{Alon-Tarsi number of $G$}. 
\end{definition}

It is clear from this definition (by using Theorem~\ref{t::CNSS} and the fact that $P_G$ is a homogeneous polynomial) that $AT(G) \geq ch(G)$ for every graph $G$. However, the converse is not true. Indeed, as was previously indicated, $ch(K_{n,n}) = (1+o(1))\log_2 n$. However, it is clear from the definition that $AT(K_{n,n}) \geq n/2$ (in fact, using different terminology, it was proved in~\cite{AT92} that $AT(K_{n,n}) = \lceil n/2 \rceil + 1$).

In theory, Theorem~\ref{t::CNSS} can be applied to many (not necessarily coloring related) combinatorial problems; however, proving that some appropriate monomial does not vanish, is often extremely hard. It is therefore convenient to give the coefficient of the monomial in question a combinatorial interpretation. The following interpretation of Alon-Tarsi choosability was proved in~\cite{AT92}. Let $G=(V,E)$ be an undirected graph with vertex set $V = \{v_1, \ldots, v_n\}$ and let $D=(V,\vec{E})$ be an orientation of $G$. Let $\vec{d} = (d_1, \ldots, d_n)$ denote the vector of out-degrees of the vertices of $D$; that is, $d_i$ is the out-degree of $v_i$ in $D$ for every $1 \leq i \leq n$. Let $EE(D)$ (respectively $EO(D)$) denote the set of all eulerian subdigraphs of $D$ (a digraph is eulerian if the out-degree of every vertex equals its in-degree; connectivity is not required --- in particular, the empty graph is considered to be an even eulerian subdigraph of $D$) with an even (respectively odd) number of edges. 

\begin{theorem} [\cite{AT92}] \label{mainAT}
Let $G=(V,E)$ be an undirected graph with vertex set $V = \{v_1, \ldots, v_n\}$. Let $D$ be an orientation of $G$ in which the out-degree of $v_i$ is $d_i$ for every $1 \leq i \leq n$. Let $f : \{1, \ldots, n\} \rightarrow \mathbb{N}$ be the function satisfying $f(i) = d_i + 1$ for every $1 \leq i \leq n$. Then, $G$ is $f-AT$ if and only if $|EE(D)| \neq |EO(D)|$.
\end{theorem}

Several successful applications of Theorem~\ref{mainAT} appear in the literature. In particular, it was used in~\cite{AT92} to prove that $AT(G) \leq 3$ whenever $G$ is a planar bipartite graph, and in~\cite{FS} to prove that $AT(G) \leq 3$ whenever $G=(V, E_1 \cup E_2)$ is a graph on $3n$ vertices such that $E_1 \cap E_2 = \emptyset$, $G_1 := (V, E_1) \cong C_n$, and $G_2 := (V, E_2)$ is a triangle factor (proving in particular a conjecture of Erd\H{o}s asserting that $G$ is $3$-colorable). 

Our first result is a generalization of Theorem~\ref{mainAT}:

\begin{theorem} \label{weightedSubgraphs}
Let $G=(V,E)$ be an undirected graph with vertex set $V = \{v_1, \ldots, v_n\}$. Let $D$ be an orientation of $G$ in which the out-degree of $v_i$ is $d_i$ for every $1 \leq i \leq n$. Let $f : \{1, \ldots, n\} \rightarrow \mathbb{N}$ be the function satisfying $f(i) = d_i + 1$ for every $1 \leq i \leq n$. For every $1 \leq i \leq n$ and every $1 \leq j \leq d_i$ let $u^i_j$ be an arbitrary real number. Then, $G$ is $f-AT$ if and only if $\sum_{A \subseteq E} (-1)^{|A|} \prod_{i=1}^n \prod_{j=1}^{d_i} \left((d^+_A(v_i) - d^-_A(v_i)) - u^i_j\right) \neq 0$.
\end{theorem}

Note that Theorem~\ref{mainAT} is a special case of Theorem~\ref{weightedSubgraphs}, obtained by setting $u^i_j = j$ for every $1 \leq i \leq n$ and every $1 \leq j \leq d_i$. Indeed, if $A \subseteq E$ does not span a eulerian subdigraph of $D$, then there must exist some $1 \leq i \leq n$ for which $d^+_A(v_i) \neq d^-_A(v_i)$. Since the sum of out-degrees is equal to the sum of in-degrees in every digraph, it follows that there must exist some $1 \leq k \leq n$ for which $1 \leq d^+_A(v_k) - d^-_A(v_k) \leq d_k$. Hence 
$$
\prod_{i=1}^n \prod_{j=1}^{d_i} \left((d^+_A(v_i) - d^-_A(v_i)) - u^i_j\right) = 0.
$$ 
If on the other hand $A \subseteq E$ does span a eulerian subdigraph of $D$, then 
$$
(-1)^{|A|} \prod_{i=1}^n \prod_{j=1}^{d_i} \left((d^+_A(v_i) - d^-_A(v_i)) - u^i_j\right) = (-1)^{|A|+m} \prod_{i=1}^n (d_i)! \ .
$$ 
It follows that 
\begin{eqnarray*}
&& \sum_{A \subseteq E} (-1)^{|A|} \prod_{i=1}^n \prod_{j=1}^{d_i} \left((d^+_A(v_i) - d^-_A(v_i)) - u^i_j\right) = \\ 
&& (-1)^m \left(|EE(D)| - |EO(D)|\right) \cdot \prod_{i=1}^n (d_i)! \ .
\end{eqnarray*}

Our proof of Theorem~\ref{weightedSubgraphs} is very different than the proof of Theorem~\ref{mainAT} given by Alon and Tarsi in~\cite{AT92}.

\bigskip

In~\cite{Alon}, Alon gave the following interpretation of $AT(G)$:

\begin{theorem} \label{signedColorings}
Let $H$ be a $d$-regular $d$-edge-colorable multi-graph and let $G$ be the line graph of $H$. Let $V(G) = \{v_1, \ldots, v_n\}$. Let $C_d$ denote the set of proper colorings of $G$ with colors $\{0,1, \ldots, d-1\}$. For $c \in C_d$, the \emph{sign} of $c$, denoted $sign(c)$, is defined to be $1$ if $P_G(c(v_1), \ldots, c(v_n)) > 0$ and $-1$ otherwise. Then, $AT(G) \leq d$ if and only if $\sum_{c \in C_d} sign(c) \neq 0$. 
\end{theorem} 

\textbf{Remark:} The definition of $sign(c)$ given in~\cite{Alon} is different than the one appearing in Theorem~\ref{signedColorings} above. However, it is not hard to see that both definitions are in fact equivalent (up to reversing all signs).

Theorem~\ref{signedColorings} was used in~\cite{EG} to prove that $AT(G) = ch(G) = \chi(G)$ whenever $G$ is the line graph of a planar $d$-regular $d$-edge-colorable multi-graph.

Our next result is a generalization of Theorem~\ref{signedColorings} (which is itself a generalization of results from~\cite{AT92, Jaeger, DES}) from certain line graphs to any graph.

\begin{theorem} \label{achFormula}
Let $G=(V,E)$ be a graph, where $V = \{1,2, \ldots, n\}$ and $|E|=m$, and let $f : V \rightarrow \mathbb{N}^{+}$ be a function satisfying $\sum_{u \in V} f(u) = m+n$. Let $C_f$ denote the set of all proper colorings $c : V \rightarrow \{0,1, \ldots, \max \{f(i)-1 : 1 \leq i \leq n\}\}$ of $G$. Then, $G$ is $f-AT$ if and only if
\begin{eqnarray*}
\sum_{c \in C_f} (-1)^{\sum_{i=1}^n c(i)} \prod_{i=1}^{n}{f(i)-1 \choose c(i)} P_G(c(1), \ldots, c(n)) \neq 0.
\end{eqnarray*}
\end{theorem}

In order to see that Theorem~\ref{achFormula} is indeed a generalization of Theorem~\ref{signedColorings}, note that if $G=(\{1, \ldots, n\}, E)$ is the line graph of a $d$-regular $d$-edge-colorable multi-graph and $c_1, c_2 : V(G) \rightarrow \{0,1, \ldots, d-1\}$ are two proper colorings, then
\begin{eqnarray*}
&& (-1)^{\sum_{i=1}^{n}c_1(i)} \prod_{i=1}^{n}{d-1 \choose c_1(i)} |P_G(c_1(1), \ldots, c_1(n))| =\\ 
&& (-1)^{\sum_{i=1}^{n}c_2(i)} \prod_{i=1}^{n}{d-1 \choose c_2(i)} |P_G(c_2(1), \ldots, c_2(n))|.
\end{eqnarray*}

Note that since $ch(G) \leq AT(G)$ for every graph $G$, it follows by Theorem~\ref{achFormula} that one can sometimes upper bound the choice number of $G$ by using only information on the proper colorings of $G$.

We will give two proofs of Theorem~\ref{achFormula}. Both proofs are very different than Alon's proof of Theorem~\ref{signedColorings} in~\cite{Alon}. The second proof will in fact entail a stronger result which we will discuss further in Section~\ref{sec::openprobs}. 

We will also prove some corollaries of Theorem~\ref{achFormula}.

An easy way to apply Theorem~\ref{achFormula} is to find lists $L(u)$ of the appropriate sizes, from which $G=(V,E)$ admits a unique coloring (that is, a unique mapping $c : V \rightarrow \bigcup_{u \in V} L(u)$ such that $f(u) \neq f(v)$ whenever $(u,v) \in E$ and $f(u) \in L(u)$ for every $u \in V$). One such example is $K_n$ with lists $\{\{0, \ldots, i-1\} : 1 \leq i \leq n\}$. This idea was explored in~\cite{AMS}. Of course one cannot expect to always find such lists. A slightly more common scenario is when the graph itself is uniquely colorable. A graph $G$ is called \emph{uniquely $k$-colorable} if there is exactly one proper vertex $k$-coloring of $G$ up to permutations of the color classes. Uniquely colorable graphs have been studied extensively and were shown to have many interesting properties (see e.g.~\cite{HHR, CG, BS}).

It is well known (see e.g.~\cite{Tru, Xu}) that a uniquely $k$-colorable graph on $n$ vertices must have at least $(k-1)n - {k \choose 2}$ edges. Moreover, there exist graphs that attain this minimum and they form a rich and interesting class (see e.g.~\cite{Tru, AMS01, Fowler, KZ, Xu}). Using Theorem~\ref{achFormula} we can easily determine their Alon-Tarsi number.

\begin{theorem} \label{theo::unique2}
Let $G=(V,E)$ be a uniquely $k$-colorable graph with $n$ vertices and $m$ edges. If $m = (k-1)n - {k \choose 2}$, then $AT(G) = k$.
\end{theorem} 

Note that if $G$ is uniquely $k$-colorable, then $\chi(G) = k$. It follows that $\chi(G) = ch(G) = AT(G)$ whenever $G=(V,E)$ is uniquely $k$-colorable and $|E| = (k-1)|V| - {k \choose 2}$. 

Another family of uniquely $k$-colorable graphs whose Alon-Tarsi number we can determine is presented in the following theorem. 

\begin{theorem} \label{theo::unique}
Let $G=(V,E)$ be a uniquely $k$-colorable graph with $n$ vertices and $m$ edges. Let $A_1, \ldots, A_k$ denote the color classes in the unique $k$-coloring of $G$. Let $r_o$ (respectively $r_e$) denote the number of parts $A_i$ of odd (respectively even) size. Let $p_o$ (respectively $p_e$) denote the number of pairs $\{(i,j) : 1 \leq i < j \leq k\}$ for which both $|A_i|$ and $|A_j|$ are odd (respectively even) and $e(A_i,A_j)$ is odd. If $m \leq \max \{(n-r_o)(k-1) + {r_o \choose 2} - p_e, (n-r_e)(k-1) + {r_e \choose 2} - p_o\}$, then $AT(G) = k$.
\end{theorem} 

While Theorem~\ref{theo::unique} might seem somewhat technical, it generalizes several known results:

\begin{corollary} [\cite{ERT}]~\label{K2*n}
$AT(K_{2*n}) = n$, where $K_{2*n}$ is the complete multi-partite graph with $n$ parts, each of size $2$.
\end{corollary}

\begin{corollary} [\cite{PW}]~\label{myPW}
Let $C_n^p$ denote the \emph{$p$th power of the $n$-cycle}, that is, the vertices of $C_n^p$ are the elements of the cyclic group $\mathbb{Z}_n$, and $(i,j) \in E(C_n^p)$ if and only if $1 \leq (i-j) \mod n \leq p$ or $1 \leq (j-i) \mod n \leq p$. If $(p+1) | n$, then $AT(C_n^p) = \chi(C_n^p) = p+1$.
\end{corollary}

\subsection{Preliminaries} \label{subsec::prelim}

Our graph-theoretic notation is standard and follows that
of~\cite{West}. In particular, we use the following.
For a graph $G$, let $V(G)$ and $E(G)$ denote its sets of vertices
and edges respectively; and let $v(G) = |V(G)|$ and $e(G) = |E(G)|$. 
For disjoint sets $A,B \subseteq V(G)$, let $E_G(A,B)$ denote the set of edges
of $G$ with one endpoint in $A$ and one endpoint in $B$, and let $e_G(A,B) = |E_G(A,B)|$.
Sometimes, if there is no risk of confusion, we discard the subscript $G$ in the above notation. Let $D=(V',E')$ be a directed graph. For an arc $(u,v) \in E'$ we say that $u$ is its \emph{tail} and $v$ is its \emph{head}. Let $A \subseteq E'$ and let $u \in V'$. We denote by $d_A^+(v)$ the number of arcs $e \in A$ such that $v$ is the tail of $e$. Similarly, we denote by $d_A^-(v)$ the number of arcs $e \in A$ such that $v$ is the head of $e$. We refer to $d_A^+(v)$ as the \emph{out-degree} of $v$ in $A$ and to $d_A^-(v)$ as the \emph{in-degree} of $v$ in $A$. Let $\mathbb{N}$ denote the set of natural numbers and let $\mathbb{N}^+ := \mathbb{N} \setminus \{0\}$. Let $\mathbb{R}$ denote the field of real numbers. Let $f, f' : X \rightarrow \mathbb{R}$ be two functions. We write $f \leq f'$ if $f(x) \leq f'(x)$ holds for every $x \in X$.   

\medskip

We mention here, without a proof, some simple observations regarding Alon-Tarsi choosability.
\begin{observation} \label{trivial}
Let $G=(V,E)$ be a multi-graph with $n$ vertices and $m$ edges, then
\begin{enumerate}
\item $P_G$ is homogeneous of degree $m$;
\item $AT(G) \geq 1 + \lceil m(G) \rceil$, where $m(G) := \max \left\{ \frac{e(H)}{v(H)} : \emptyset \neq H \subseteq G \right\}$ is the maximum density of $G$;
\item If $C_1, \ldots, C_r$ are the connected components of $G$, then $P_G = \prod_{i=1}^r P_{C_i}$ and, in particular, $AT(G) = \max \{AT(C_i) : 1 \leq i \leq r\}$;
\item If $H$ is a subgraph of $G$, then $AT(H) \leq AT(G)$ (note that, unlike with colorability or choosability, a strict inequality can occur even when $H$ is obtained from $G$ by replacing some parallel edges with a single edge).
\item Let $f, f' : V \rightarrow \mathbb{N}^+$ be two functions satisfying $f' \geq f$. If $G$ is $f-AT$, then it is also $f'-AT$. 
\item $AT(G) \leq col(G) + 1$, where $col(G) = \max \{\delta(H) : H \subseteq G\}$ is the coloring number of $G$.
\end{enumerate}
\end{observation}

\textbf{Remark:} In light of part $3.$ of Observation~\ref{trivial}, henceforth we will consider only connected multi-graphs.

\noindent The rest of the paper is organized as follows: in
Section~\ref{sec::proofs} we prove
Theorems~\ref{weightedSubgraphs} and~\ref{achFormula}. In Section~\ref{sec::unique} we prove Theorems~\ref{theo::unique2} and~\ref{theo::unique}, as well as Corollaries~\ref{K2*n} and~\ref{myPW}. Finally, in Section~\ref{sec::openprobs} we discuss some generalizations and present some open problems.

\section{Generalizing the polynomial method} \label{sec::proofs}

\noindent \textbf{Proof of Theorem~\ref{weightedSubgraphs}:}

By Theorem~\ref{t::CNSS} it suffices to prove that the coefficient of $\prod_{i=1}^n x_i^{d_i}$ in the expansion of $P_G(x_1, \ldots, x_n)$ is non-zero if and only if 
\begin{eqnarray*}
\sum_{A \subseteq E} (-1)^{|A|} \prod_{i=1}^n \prod_{j=1}^{d_i} \left((d^+_A(v_i) - d^-_A(v_i)) - u^i_j\right) \neq 0. 
\end{eqnarray*}
Let $E = \{e_1, \ldots, e_m\}$. Let $A = A_D = (a_{i,j})_{m \times n}$ denote the oriented incidence matrix of $D$, that is, for every $1 \leq i \leq m$ and $1 \leq j \leq n$ we have $a_{i,j} = 1$ if $j$ is the tail of $e_i$, $a_{i,j} = -1$ if $j$ is the head of $e_i$, and $a_{i,j} = 0$ otherwise. It is easy to see that $\left|P_G(x_1, \ldots, x_n)\right| = \left|\prod_{i=1}^m \left(\sum_{j=1}^n a_{i,j} x_j\right)\right|$. Let $B = (b_{i,j})_{m \times m}$ be the matrix obtained from $A$ by repeating $A$'s $j$th column exactly $d_j$ times for every $1 \leq j \leq n$, that is, for every $1 \leq p \leq n$ and every $1 \leq q \leq d_p$, the $(\sum_{i=1}^{p-1} d_i + q)$th column of $B$ is equal to the $p$th column of $A$. It is not hard to see (see e.g. Claim 1 in~\cite{AT89}) that the coefficient of $\prod_{i=1}^n x_i^{d_i}$ in the expansion of $P_G(x_1, \ldots, x_n)$ does not vanish if and only if $Per(B) \neq 0$. Finally, let $C = (c_{i,j})_{(m+1) \times (m+1)}$ be the matrix satisfying $c_{i,j} = b_{i,j}$ for every $1 \leq i,j \leq m$, $c_{i,m+1} = 0$ for every $1 \leq i \leq m$, $c_{m+1,m+1} = 1$, and $c_{m+1,j} = - u^p_q$ whenever $1 \leq j = \sum_{i=1}^{p-1} d_i + q \leq m$ for some $1 \leq p \leq n$ and $1 \leq q \leq d_p$. Expanding the permanent of $C$ using its $(m+1)$st column, it is clear that $Per(C) = Per(B)$. It follows that the coefficient of $\prod_{i=1}^n x_i^{d_i}$ in the expansion of $P_G(x_1, \ldots, x_n)$ does not vanish if and only if $Per(C) \neq 0$. In order to compute the permanent of $C$ we will use Ryser's formula~\cite{Ryser} and the simple fact that a matrix and its transpose have the same permanent:
\begin{eqnarray*} 
(-1)^{m+1} Per(C) &=& (-1)^{m+1} Per(C^T)\\ 
&=& \sum_{S \subseteq \{1, \ldots, m+1\}} (-1)^{|S|} \prod_{k=1}^{m+1} \left(\sum_{r \in S} c_{r,k}\right)\\
&=& \sum_{S \subseteq \{1, \ldots, m+1\} \atop m+1 \in S} (-1)^{|S|} \prod_{k=1}^m \left(\sum_{r \in S} c_{r,k}\right) \cdot 1\\ 
&+& \sum_{S \subseteq \{1, \ldots, m\}} (-1)^{|S|} \prod_{k=1}^m \left(\sum_{r \in S} c_{r,k}\right) \cdot 0\\
&=& \sum_{A \subseteq E} (-1)^{|A|} \prod_{i=1}^n \prod_{j=1}^{d_i} \left((d^+_A(v_i) - d^-_A(v_i)) - u^i_j\right).
\end{eqnarray*}
{\hfill $\Box$ \medskip\\}

\noindent \textbf{First proof of Theorem~\ref{achFormula}:}

This proof is similar in spirit to the proof of Theorem~\ref{weightedSubgraphs}. Let $D$ be the orientation of $G$ with no decreasing arcs, that is, $(i,j) \in E$ is directed from $i$ to $j$ in $D$ if and only if $i < j$. Let $A$ denote the oriented incidence matrix of $D$. Let $B = (b_{i,j})$ denote the $m \times m$ matrix which is obtained from $A$ by repeating $A$'s $j$th column exactly $f(j)-1$ times for every $1 \leq j \leq n$. By Ryser's formula we have:
\begin{eqnarray*}
(-1)^m Per(B) &=& \sum_{S \subseteq \{1, \ldots, m\}} (-1)^{|S|} \prod_{r=1}^m \left(\sum_{k \in S} b_{r,k}\right).
\end{eqnarray*}
Consider an arbitrary fixed set $S \subseteq \{1, \ldots, m\}$. Each $i \in S$ corresponds to some column of $B$, which in turn corresponds to a column of $A$ and thus to a vertex of $G$. Let $c_S : V \rightarrow \mathbb{N}$ be the mapping which assigns to every vertex $j \in V$ the number of elements $i \in S$ such that the $i$th column of $B$ corresponds to $j$. Note that $0 \leq c_S(j) \leq f(j)-1$ holds for every $j \in V$ and for every $S \subseteq \{1, \ldots, m\}$. Let $E = \{e_1, \ldots, e_m\}$, and for every $1 \leq r \leq m$ let $i_r$ denote the tail of $e_r$ (in $D$) and let $j_r$ denote the head of $e_r$. Let $S \subseteq \{1, \ldots, m\}$ be such that $c_S$ is not a proper coloring of $G$. Let $e_{\ell} \in E$ be some edge such that $c_S(i_{\ell}) = c_S(j_{\ell})$. It follows that $\sum_{k \in S} b_{\ell,k} = c_S(i_{\ell}) - c_S(j_{\ell}) = 0$ and thus $(-1)^{|S|} \prod_{r=1}^m \left(\sum_{k \in S} b_{r,k}\right) = 0$. If on the other hand $S \subseteq \{1, \ldots, m\}$ is such that $c_S$ is a proper coloring of $G$, then 
\begin{eqnarray*}
(-1)^{|S|} \prod_{r=1}^m \left(\sum_{k \in S} b_{r,k}\right) &=& (-1)^{|S|} \prod_{r=1}^m (c_S(i_r) - c_S(j_r))\\
&=& (-1)^{\sum_{i=1}^n c_S(i)} P_G(c_S(1), \ldots, c_S(n)).
\end{eqnarray*}
Since every coloring $c_S$ is obtained from exactly $\prod_{i=1}^n {f(i)-1 \choose c_S(i)}$ different sets $S \subseteq \{1, \ldots, m\}$ (as every choice of exactly $c_i$ columns of $B$ which corresponds to $i \in V$, for every $1 \leq i \leq n$, yields the same coloring), it follows that 
\begin{eqnarray*}
&& \sum_{S \subseteq \{1, \ldots, m\}} (-1)^{|S|} \prod_{r=1}^m \left(\sum_{k \in S} b_{r,k}\right) = \\ 
&& \sum_{c} (-1)^{\sum_{i=1}^n c(i)} \prod_{i=1}^n {f(i)-1 \choose c(i)} P_G(c(1), \ldots, c(n)),
\end{eqnarray*}
where the second sum extends over all proper colorings $c$ of $G$ in which $0 \leq c(i) \leq f(i)-1$ holds for every $1 \leq i \leq n$. Finally, since ${f(i)-1 \choose c(i)} = 0$ whenever $c(i) > f(i)-1$, it follows that
\begin{eqnarray*}
(-1)^m Per(B) = \sum_{c \in C_f} (-1)^{\sum_{i=1}^n c(i)} \prod_{i=1}^{n}{f(i)-1 \choose c(i)} P_G(c(1), \ldots, c(n)).
\end{eqnarray*}
As in the proof of Theorem~\ref{weightedSubgraphs}, it follows that the coefficient of $\prod_{i=1}^n x_i^{f(i)-1}$ in the expansion of $P_G(x_1, \ldots, x_n)$ does not vanish if and only if 
\begin{eqnarray*}
\sum_{c \in C_f} (-1)^{\sum_{i=1}^n c(i)} \prod_{i=1}^{n}{f(i)-1 \choose c(i)} P_G(c(1), \ldots, c(n)) \neq 0. 
\end{eqnarray*}
The result now follows by Theorem~\ref{t::CNSS}.
{\hfill $\Box$ \medskip\\}

\noindent \textbf{Second proof of Theorem~\ref{achFormula}:}

We will prove the following more general theorem:

\begin{theorem} \label{genListVersion}
Let $P \in \mathbb{R}[x_1, \ldots, x_n]$ be a polynomial of degree $m$. Let $f : \{1, \ldots, n\} \rightarrow \mathbb{N}^+$ be a function satisfying $\sum_{i=1}^n f(i) = m+n$. Let $NZ(P) := \{(z_1, \ldots, z_n) \in \{0, \ldots, f(1)-1\} \times \ldots \times \{0, \ldots, f(n)-1\} : P(z_1, \ldots, z_n) \neq 0\}$ denote the set of \emph{non-zeros} of $P$ in $\{0, \ldots, f(1)-1\} \times \ldots \times \{0, \ldots, f(n)-1\}$. If
\begin{eqnarray*}
\sum_{(z_1, \ldots, z_n) \in NZ(P)} (-1)^{\sum_{i=1}^n z_i} \prod_{i=1}^n \binom{f(i)-1}{z_i} P(z_1, \ldots, z_n) \neq 0,
\end{eqnarray*}
then, for every $S_1, \ldots, S_n \subseteq \mathbb{R}$, where $|S_i| \geq f(i)$ for every $1 \leq i \leq n$, there exists $(s_1, \ldots, s_n) \in S_1 \times \ldots \times S_n$ such that $P(s_1, \ldots, s_n) \neq 0$.
\end{theorem}

Theorem~\ref{achFormula} is the special case of Theorem~\ref{genListVersion} obtained by taking $P = P_G$ (the ``only if'' part of Theorem~\ref{achFormula} follows directly from the definition of Alon-Tarsi choosability). Indeed, $NZ(P_G)$ are exactly the proper colorings of $G$ (observe that ${f(i)-1 \choose c(i)} = 0$ whenever $c(i) > f(i)-1$). Other special cases of Theorem~\ref{genListVersion} will be discussed in Section~\ref{sec::openprobs}.

Theorem~\ref{genListVersion} is a fairly straightforward corollary of Theorem~\ref{t::CNSS} and of the following result:

\begin{lemma} [Scheim~\cite{DES}] \label{lem::scheim}
\label{scheim} If $P(x_{1},x_{2},\ldots,x_{n}) \in
\mathbb{R}[x_{1},x_{2},\ldots,x_{n}]$ is of degree $\leq s_{1} +
s_{2} + \ldots + s_{n}$, where $n$ is a positive integer and
$s_{1},s_{2},\ldots,s_{n}$ are non-negative integers, then
\begin{eqnarray*}
&&\left(\frac{\partial}{\partial x_{1}}\right)^{s_{1}}\left(\frac{\partial}{\partial
x_{2}}\right)^{s_{2}} \ldots \left(\frac{\partial}{\partial
x_{n}}\right)^{s_{n}}P(x_{1},x_{2},\ldots,x_{n})\\ &=& 
\sum_{x_{1}=0}^{s_{1}} \ldots
\sum_{x_{n}=0}^{s_{n}}(-1)^{s_{1}+x_{1}}{s_{1} \choose
x_{1}}\ldots(-1)^{s_{n}+x_{n}}{s_{n} \choose
x_{n}}P(x_{1},x_{2},\ldots,x_{n}).
\end{eqnarray*}
\end{lemma}

\noindent \textbf{Proof of Theorem~\ref{genListVersion}:}
By Theorem~\ref{t::CNSS} it suffices to prove that $c_f \neq 0$, where $c_f$ is the coefficient of $\prod_{i=1}^n x_i^{f(i)-1}$ in the expansion of $P(x_1, \ldots, x_n)$. However, by Lemma~\ref{scheim} we have
\begin{eqnarray*}
c_f \prod_{i=1}^n (f(i)-1)! = \sum_{x_1=0}^{f(1)-1} \ldots \sum_{x_n=0}^{f(n)-1}(-1)^m (-1)^{\sum_{i=1}^n x_i} \prod_{i=1}^n {f(i)-1 \choose x_i} P(x_1, x_2, \ldots, x_n).
\end{eqnarray*}
Thus, in order to prove that $c_f \neq 0$, it suffices to show that
\begin{eqnarray*}
\sum_{x_1=0}^{f(1)-1} \ldots \sum_{x_n=0}^{f(n)-1} (-1)^{\sum_{i=1}^n x_i} \prod_{i=1}^n {f(i)-1 \choose x_i} P(x_1, x_2, \ldots, x_n) \neq 0.
\end{eqnarray*}
By the definition of $NZ(P)$, we have $P(x_1, \ldots, x_n) = 0$ whenever $(x_1, \ldots, x_n) \in \{0, \ldots, f(1)-1\} \times \ldots \times \{0, \ldots, f(n)-1\} \setminus NZ(P)$. It follows that this last expression equals
\begin{eqnarray*}
\sum_{(x_1, \ldots, x_n) \in NZ(P)} (-1)^{\sum_{i=1}^n x_i} \prod_{i=1}^n \binom{f(i)-1}{x_i} P(x_1, x_2, \ldots, x_n).
\end{eqnarray*}
This concludes the proof of the theorem as this last expression was assumed to be non-zero.
{\hfill $\Box$ \medskip\\}

\section{Uniquely colorable graphs} \label{sec::unique}

\noindent \textbf{Proof of Theorem~\ref{theo::unique2}:}
Clearly, $AT(G) \geq \chi(G) = k$. Let $A_0, A_1, \ldots, A_{k-1}$ be the color classes of the unique $k$-coloring of $G$. Arbitrarily choose a vertex $v_i \in A_i$ for every $0 \leq i \leq k-1$. Assign the list $L(v_i) = \{0, \ldots, i\}$ to $v_i$ for every $0 \leq i \leq k-1$ and the list $L(u) = \{0, \ldots, k-1\}$ to every vertex $u \in V \setminus \{v_0, \ldots, v_{k-1}\}$. It is easy to see that $\sum_{u \in V} |L(u)| = m+n$ and that the only proper coloring of $G$ from these lists is the one assigning the color $i$ to every vertex of $A_i$ for $0 \leq i \leq k-1$. The theorem now readily follows from Theorem~\ref{achFormula}.
{\hfill $\Box$ \medskip\\}

\textbf{Remark:} One class to which Theorem~\ref{theo::unique2} applies, is that of uniquely $4$-colorable planar graphs. However, for such graphs, the inequality $AT(G) \leq 4$ follows immediately since they are $3$-degenerate~\cite{Fowler}. On the other hand, as was partly indicated in the Introduction, for every $k \geq 3$ there are infinitely many uniquely $k$-colorable graphs with the minimum possible number of edges which are \emph{not} $(k-1)$-degenerate (see e.g.~\cite{Tru, AMS01, Fowler, KZ, Xu}).

\bigskip

\noindent \textbf{Proof of Theorem~\ref{theo::unique}:}
Clearly, $AT(G) \geq \chi(G) = k$. Assume without loss of generality that $(n-r_o)(k-1) + {r_o \choose 2} - p_e \geq (n-r_e)(k-1) + {r_e \choose 2} - p_o$, and let $r = r_o$ and $p = p_e$. Assume without loss of generality that $A_1, \ldots, A_r$ are the odd parts. For every $r+1 \leq i < j \leq k$ for which $e(A_i,A_j)$ is odd, add an arbitrary edge (allowing parallel edges) between $A_i$ and $A_j$; denote the resulting multi-graph by $G'$. By our assumption we will add exactly $p$ such edges; let $m' = m+p$ denote the number of edges in $G'$. For every $1 \leq i \leq r$, let $u_i$ be an arbitrary vertex of $A_i$. Assign the list $L(u_i) = \{0, \ldots, i-1\}$ to $u_i$ for every $1 \leq i \leq r$ and the list $L(u) = \{0, \ldots, k-1\}$ to every vertex $u \in V \setminus \{u_1, \ldots, u_r\}$. It is easy to see that $\sum_{u \in V} |L(u)| \geq m' + n$. Add $s := \sum_{u \in V} |L(u)| - (m' + n)$ arbitrary edges between $A_1$ and $A_2$ (allowing parallel edges); denote the resulting multi-graph by $G''$. Note that $\sum_{u \in V} |L(u)| = e(G'') + n$. Moreover, for any $r+1 \leq i < j \leq k$ there is an even number of edges between $A_i$ and $A_j$. Indeed, this holds in $G'$ by construction. If $r \geq 1$, then this holds trivially in $G''$. If $r=0$, then $n$ is even, $m'$ is even, and $\sum_{u \in V} |L(u)| = kn$ is even. It follows that, in this case, $s$ is even as well. It is evident that $c : V \rightarrow \{0, \ldots, k-1\}$ assigning color $i-1$ to every vertex of $A_i$ for every $1 \leq i \leq k$, is a proper coloring of $G''$ from the lists $\{L(u) : u \in V\}$. Moreover, it is clear that any proper coloring $c'$ of $G''$ from these lists is sign preserving, that is
\begin{eqnarray*}
(-1)^{\sum_{i=1}^n c'(i)} sign(P_G(c'(1), \ldots, c'(n))) = (-1)^{\sum_{i=1}^n c(i)} sign(P_G(c(1), \ldots, c(n))). 
\end{eqnarray*}
It follows by Theorem~\ref{achFormula} and Observation~\ref{trivial} that $AT(G) \leq AT(G') \leq AT(G'') \leq k$.
{\hfill $\Box$ \medskip\\}

\noindent \textbf{Proof of Corollary~\ref{K2*n}:}
It is clear that $K_{2*n}$ is uniquely $n$-colorable, that all color classes are of even size, and that the number of edges of $K_{2*n}$ between any two color classes is even. Since $e(K_{2*n}) = 2n(n-1)$, it follows from Theorem~\ref{theo::unique} that $AT(K_{2*n}) \leq n$.
{\hfill $\Box$ \medskip\\}

\noindent \textbf{Proof of Corollary~\ref{myPW}:}
It is easy to see that if $(p+1) | n$, then $C_n^p$ is uniquely $(p+1)$-colorable and that all color classes in this unique coloring are of the same size $n/(p+1)$. Moreover, the number of edges between any two color classes is even, as every vertex of $C_n^p$ has exactly $2$ neighbors from every color class (other than its own). The result now follows from Theorem~\ref{theo::unique}, as clearly $e(C_n^p) = pn$.
{\hfill $\Box$ \medskip\\}

\textbf{Remark:} Prowse and Woodall~\cite{PW} have determined $AT(C_n^p)$ for every $n$ and $p$ and not only when $(p+1) | n$. However, their proof of this special case is also quite involved.

\section{Concluding remarks and open problems} \label{sec::openprobs}

\textbf{More special cases of Theorem~\ref{genListVersion}.} This theorem demonstrates how the use of the Combinatorial Nullstellensatz (Theorem~\ref{t::CNSS}) together with Scheim's Lemma (Lemma~\ref{lem::scheim}) can provide a general framework for transforming certain labeling problems into their ``list coloring version''. That is, under appropriate circumstances, one can substitute a labeling of some object with elements from $\{0,1, \ldots, d_1\} \times \ldots \times \{0,1, \ldots, d_n\}$ that satisfies some property ${\mathcal P}$, with a labeling of the same object which also satisfies ${\mathcal P}$ but uses the elements of $S_1 \times \ldots \times S_n$, where, for every $1 \leq i \leq n$, $S_i$ is an arbitrary set of size $d_i + 1$. This idea was already used in~\cite{Hefetz}. One can think of many such examples. We mention two of them here:

\textbf{Alon-Tarsi choosability of hypergraphs:} The notion of list coloring easily extends to hypergraphs; however, very little is known about choice numbers of hypergraphs.  

Ramamurthi and West~\cite{RWest} have extended the Alon-Tarsi polynomial method to $k$-uniform hypergraphs, where $k$ is prime. Their definition of the \emph{hypergraph polynomial} is as follows. Let $k$ be a prime, let $\omega$ be a primitive $k$th root of unity, and let ${\mathcal H} = (V,E)$ be a $k$-uniform hypergraph, where $V = \{1,2, \ldots, n\}$. With every vertex $i \in V$ associate a variable $x_i$. For an edge $e = \{i_0, \ldots, i_{k-1}\}$ of ${\mathcal H}$, where $i_0 < \ldots < i_{k-1}$, define $P_e = \sum_{j=0}^{k-1}\omega^j x_{i_j}$. The hypergraph polynomial of ${\mathcal H}$ is defined by $P_{\mathcal H}(x_1, \ldots, x_n) = \prod_{e \in E} P_e$. As was indicated in~\cite{RWest}, ${\mathcal H}$ is $r$-colorable if and only if there exists an $n$-tuple $(a_1, \ldots, a_n) \in \{0, \ldots, r-1\}^n$ such that $P_{\mathcal H}(a_1, \ldots, a_n) \neq 0$. Using this definition, one can easily generalize Definition~\ref{algChoosable} to $k$-uniform hypergraphs, where $k$ is prime. However, as observed in~\cite{pei}, in order to upper bound $AT({\mathcal H})$, it suffices to associate ${\mathcal H}$ with any polynomial $Q_{\mathcal H}$ such that, if $Q_{\mathcal H}(c(1), \ldots, c(n)) \neq 0$, then $c$ is a proper coloring of ${\mathcal H}$. Using this observation we can obtain a theorem which applies to any (not necessarily uniform) hypergraph. 

\begin{theorem} \label{achHypergraph}
Let ${\mathcal H}=(V,E)$ be any hypergraph, where $V = \{1,2, \ldots, n\}$ and $|E|=m$. Let $f : V \rightarrow \mathbb{N}^+$ be a function satisfying $\sum_{u \in V} f(u) = m+n$. Let $Q = Q_{{\mathcal H}, f}(x_1, \ldots, x_n)$ be any polynomial of degree $m$ such that $Q(c(1), \ldots, c(n)) = 0$ whenever $(c(1), \ldots, c(n)) \in \{0, \ldots, f(1)-1\} \times \ldots \times \{0, \ldots, f(n)-1\}$ is \emph{not} a proper coloring of ${\mathcal H}$. Let $NZ(Q)$ denote the set of \emph{non-zeros} of $Q$ in $\{0, \ldots, f(1)-1\} \times \ldots \times \{0, \ldots, f(n)-1\}$. If
\begin{eqnarray*}
\sum_{(c_1, \ldots, c_n) \in NZ(Q)} (-1)^{\sum_{i=1}^n c_i} \prod_{i=1}^{n}{f(i)-1 \choose c_i} Q(c_1, \ldots, c_n) \neq 0,
\end{eqnarray*}
then ${\mathcal H}$ is $f-AT$.
\end{theorem} 

We demonstrate the convenience of choosing a polynomial, by proving that $AT(F_7) \leq 3$, where $F_7$ is the Fano plane (it is instructive to compare this proof with the one given in~\cite{RWest}). Let 
$$
F_7 = (\{1,2, \ldots, 7\}, \{\{1,2,4\}, \{1,3,6\}, \{1,5,7\}, \{2,3,5\}, \{3,4,7\}, \{2,6,7\}, \{4,5,6\}\}),
$$ 
and let 
\begin{eqnarray*}
Q_{F_7} &=& (x_2 - x_4)(x_1 + x_3 - 2x_6)(x_1 + x_5 - 2x_7)(x_2 + x_3 - 2x_5)(x_3 - x_4)\\
&& (x_2 + x_6 - 2x_7)(x_4 + x_5 - 2x_6).
\end{eqnarray*}

It is easy to see that if $Q_{F_7}(c(1), \ldots, c(7)) \neq 0$, then $c$ is a proper coloring of $F_7$ (though the converse does not hold). Define $f : \{1,2, \ldots, 7\} \rightarrow \mathbb{N}^+$ by $f(1)=3, f(3)=1$ and $f(i)=2$ for every $i \in \{1,2, \ldots, 7\} \setminus \{1,3\}$. Note that $\sum_{i=1}^7 f(i) = 14 = e(F_7) + v(F_7)$. It is easy to see that the only proper (with respect to $Q_{F_7}$ and $f$) coloring of $F_7$ is the one defined by $c(1)=2, c(2) = c(3) = c(6) = 0$, and $c(4) = c(5) = c(7) = 1$. It follows by Theorem~\ref{achHypergraph} that $F_7$ is $f-AT$; in particular $AT(F_7) \leq 3$.

Ramamurthi and West~\cite{RWest} have noted that it would be interesting to use the ``extended Alon-Tarsi polynomial method'' to determine the choice number of an infinite family of hypergraphs for which the choice number is not yet known. Theorem~\ref{achHypergraph} might prove useful in this respect.

\textbf{$T$-list coloring:} Let $\emptyset \neq T \subseteq \mathbb{N}$ be a set. A \emph{$T$-coloring} of a graph $G=(V,E)$ with $k$ colors is a mapping $c : V \rightarrow \{0,1, \ldots, k-1\}$ such that $|c(u) - c(v)| \notin T$ whenever $(u,v) \in E$. Given a graph $G=(V,E)$, a set $T$ as above and a family ${\mathcal L} = \{L(u) \subseteq \mathbb{N} : u \in V\}$, an \emph{${\mathcal L}$-list $T$-coloring} of $G$ is a $T$-coloring $c$ of $G$ such that $c(u) \in L(u)$ for every $u \in V$. A graph is said to be \emph{$T$-$k$-choosable} if it admits an ${\mathcal L}$-list $T$-coloring whenever the family ${\mathcal L}$ contains only sets of size at least $k$. Note that when $T =\{0\}$, the notion of $T$-coloring reduces to that of coloring and the notion of $T$-list coloring reduces to that of list coloring. $T$-colorings and $T$-list colorings were extensively studied in~\cite{Tes1, Tes2}. In~\cite{AZ}, a simple variation of Theorem~\ref{mainAT} was used to prove (among other things) that even cycles are $T$-$2|T|$-choosable for every set $T$ which contains $0$. We will show how this result can be obtained using Theorem~\ref{genListVersion} instead. Let $T \subseteq \mathbb{N}$ be an arbitrary set of size $\ell$, where $0 \in T$. Let $V(C_{2n}) = \{1, \ldots, 2n\}$. Associate a variable $x_i$ with every $1 \leq i \leq 2n$. Let 
\begin{eqnarray*}    
Q_{C_{2n},T}(x_1, \ldots, x_{2n}) = \prod_{(i,j) \in E(C_{2n}) \atop i < j} (x_i - x_j) \prod_{t \in T \setminus \{0\}} (x_i - x_j - t)(x_i - x_j + t).
\end{eqnarray*}
It is evident that if for every family $\{S_i \subseteq \mathbb{N} : 1 \leq i \leq 2n\}$, each of size $2\ell$, there exists a vector $(a_1, \ldots, a_{2n}) \in S_1 \times \ldots \times S_{2n}$ such that $Q_{C_{2n},T}(a_1, \ldots, a_{2n}) \neq 0$, then $C_{2n}$ is $T$-$2\ell$-choosable. By Theorem~\ref{t::CNSS}, it suffices to prove that the coefficient of $\prod_{i=1}^{2n} x_i^{2\ell-1}$ in the expansion of $Q_{C_{2n},T}$ does not vanish. Note that for any set $\{s_1, \ldots\, s_{2\ell-1} \subseteq \mathbb{N}\}$, the coefficient of $\prod_{i=1}^{2n} x_i^{2\ell-1}$ in the expansion of 
\begin{eqnarray*}    
P := \prod_{(i,j) \in E(C_{2n}) \atop i < j} \prod_{k=1}^{2\ell-1} (x_i - x_j - s_k).
\end{eqnarray*}
is the same as that coefficient in the expansion of $Q_{C_{2n},T}$. Using Theorem~\ref{genListVersion} with the polynomial $P_S$, obtained by setting $s_i = i$ for every $1 \leq i \leq 2\ell-1$, we conclude that it suffices to prove that
\begin{eqnarray*}
\sum_{(z_1, \ldots, z_{2n}) \in NZ(P_S)} (-1)^{\sum_{i=1}^{2n} z_i} \prod_{i=1}^{2n} \binom{2\ell-1}{z_i} P_S(z_1, \ldots, z_{2n}) \neq 0,
\end{eqnarray*} 
where it is understood that $NZ(P_S) \subseteq \{0,1, \ldots, 2\ell-1\}^{2n}$. It is easy to see that $NZ(P_S) = \{(a, \ldots, a) : 0 \leq a \leq 2\ell-1\}$. Hence
\begin{eqnarray*}
&&\sum_{(z_1, \ldots, z_{2n}) \in NZ(P_S)} (-1)^{\sum_{i=1}^{2n} z_i} \prod_{i=1}^{2n} \binom{2\ell-1}{z_i} P_S(z_1, \ldots, z_{2n}) \\ 
&=& \sum_{a=1}^{2\ell-1} (-1)^{2na} \prod_{i=1}^{2n} \binom{2\ell-1}{a} \left[-(2\ell-1)! \right]^{2n} \\
&\neq& 0,
\end{eqnarray*}
and the claim follows.

\bigskip 

\textbf{Alon-Tarsi choosability vs. regular choosability.} As was noted in the Introduction, there are graphs $G$ for which $AT(G) > ch(G)$. Moreover, $AT(G)$ has certain distinct features, such as being sensitive to the addition of parallel edges. It is therefore interesting to study the Alon-Tarsi number of a graph as a separate graph invariant. For example, analogously to a classical result from~\cite{ERT}, it is easy to characterize the graphs $G$ which satisfy $AT(G) \leq 2$ (these are exactly the graphs whose core is an even cycle or a single vertex). On the other hand, it would be interesting to prove \emph{Alon-Tarsi strengthenings} of known results regarding choosability. Corollary~\ref{K2*n} is one such example. In~\cite{HKS}, an Alon-Tarsi strengthening of Brooks' Theorem was proved, that is, $AT(G) \leq \Delta(G)$ holds for any graph $G$ which is not an odd cycle or a complete graph. Since $ch(G) \leq AT(G) \leq col(G) + 1$ holds for every graph $G$ (see part $6.$ of Observation~\ref{trivial}), it follows that $5 \leq \max \{AT(G) : G \textrm{ is planar} \} \leq 6$. It would be interesting to decide which of the two bounds is the correct answer.  

\textbf{Scheim's Lemma.} Lemma~\ref{lem::scheim} is also useful in its own right. For example it easily entails Ryser's formula for computing the permanent of a matrix~\cite{Ryser}. Indeed, let $A = (a_{i,j})$ be an $n \times n$ matrix. It is well known (see e.g. Claim 1 in~\cite{AT89}) that the permanent of $A$ is exactly the coefficient of $\prod_{i=1}^n x_i$ in the expansion of $P_A(x_1, \ldots, x_n) := \prod_{i=1}^n (\sum_{j=1}^n a_{i,j} x_j)$. However, by Scheim's Lemma, this coefficient is precisely
\begin{eqnarray*}
&&\left(\frac{\partial}{\partial x_{1}}\right)\left(\frac{\partial}{\partial
x_{2}}\right) \ldots \left(\frac{\partial}{\partial
x_{n}}\right) P_A(x_1, x_2, \ldots, x_n)\\ 
&=& \sum_{x_{1}=0}^{1} \ldots
\sum_{x_{n}=0}^{1} (-1)^{1 + x_1} \binom{1}{x_1} \ldots (-1)^{1 + x_n} 
\binom{1}{x_n} P_A(x_1, x_2, \ldots, x_n)\\
&=& (-1)^n \sum_{S \subseteq \{1,2, \ldots n\}} (-1)^{|S|} \prod_{i=1}^n \left(\sum_{j \in S} a_{i,j}\right).
\end{eqnarray*}  

It would be interesting to obtain additional corollaries of this lemma.


\begin{thebibliography}{alpha}

\bibitem{AMS01}
S. Akbari, V. S. Mirrokni and B. S. Sadjad, $K_r$-free uniquely vertex colorable graphs with minimum possible edges, Journal of Combinatorial Theory Ser. B. 82 (2001), 316--318.

\bibitem{AMS}
S. Akbari, V. S. Mirrokni and B. S. Sadjad, A relation between choosability and uniquely list colorability, Journal of Combinatorial Theory Ser. B. 96 (2006), 577--583.

\bibitem{Alon}
N. Alon, Restricted colorings of graphs, in Surveys in Combinatorics 1993, London
Math. Soc. Lecture Notes Series 187 (K. Walker, ed.), Cambridge Univ. Press, 1993,
1-–33.

\bibitem{CNSS}
N. Alon, Combinatorial Nullstellensatz, Combinatorics, Probability and Computing 8 (1999) no. 1--2, 7--29.

\bibitem{AT89}
N. Alon and M. Tarsi, A nowhere-zero point in linear mappings, Combinatorica 9 (1989), 393--395. 

\bibitem{AT92}
N. Alon and M. Tarsi, Colorings and orientations of graphs, Combinatorica 12 (1992), 125--134.

\bibitem{AT97}
N. Alon and M. Tarsi, A note on graph colorings and graph polynomials, Journal of Combinatorial Theory Ser. B. 70 (1997), 197--201.

\bibitem{AZ}
N. Alon and A. Zaks, $T$-choosability in graphs, Discrete Applied Math. 82 (1998), 1--13.

\bibitem{BS}
B. Bollob\'as and N. Sauer, Uniquely colourable graphs with large girth, Canadian Journal of Mathematics 28 (1976), 1340--1344.

\bibitem{CG}
G. Chartrand and D. Geller, Uniquely colorable planar graphs, Journal of Combinatorial Theory, 6 (1969), 271--278.

\bibitem{EG}
M. Ellingham and L.A. Goddyn, List edge colorings of some regular planar multigraphs, Combinatorica 16 (1996), 343--352.

\bibitem{ERT}
P. Erd\H{o}s, A. L. Rubin and H. Taylor, Choosability in graphs, Congressus Numerantium
26 (1979), 125–-157.

\bibitem{FS}
H. Fleischner and M. Stiebitz, A solution to a colouring problem of P. Erd\H{o}s, Discrete Math. 101 (1992), 39--48.

\bibitem{Fowler}
T. G. Fowler, Unique coloring of planar graphs, manuscript.

\bibitem{Galvin}
F. Galvin, The list chromatic index of a bipartite multigraph, Journal of Combinatorial Theory Ser. B. 63 (1995), 153--158.

\bibitem{HHR}
F. Harary, S. T. Hedetniemi and R. W. Robinson, Uniquely colorable graphs, Journal of Combinatorial Theory, 6 (1969), 264--270.

\bibitem{Hefetz}
D. Hefetz, Anti-magic graphs via the Combinatorial Nullstellensatz, Journal of Graph Theory 50(4) (2005), 263--272.

\bibitem{HKS}
J. Hladk\'y, D. Kr\'al and U. Schauz, Algebraic proof of Brooks' Theorem, manuscript.

\bibitem{Jaeger}
F. Jaeger, On the Penrose number of cubic diagrams, Discrete Math. 74 (1989), 85--97.

\bibitem{Kahn}
J. Kahn, Asymptotically good list-colorings, Journal of Combinatorial Theory Ser. A. 73
(1996), 1--59.

\bibitem{KZ}
L. Kang and M. Zhao, Construction of uniquely vertex $k$-colorable graphs with minimum possible size, Journal of Shanghai University 11(5) (2007), 449--450.

\bibitem{LiLi}
S. Y. R. Li and W. C. W. Li, Independence numbers of graphs and generators of ideals,
Combinatorica 1 (1981), 55--61.

\bibitem{Lovasz}
L. Lov\'asz, Bounding the independence number of a graph, in: (A. Bachem, M. Gr\"otschel and B.
Korte, eds.), Bonn Workshop on Combinatorial Optimization, Annals of Discrete Mathematics
16 (1982), North Holland, Amsterdam.

\bibitem{Toft}
T. Jensen and B. Toft, \textbf{Graph Coloring Problems}, Wiley, New York, 1995.

\bibitem{pei}
M. Pei, List colouring hypergraphs and extremal results for acyclic graphs, manuscript.

\bibitem{Petersen}
J. Petersen, Die Theorie der regularen Graphs, Acta Math. 15 (1891), 193--220.

\bibitem{PW}
A. Prowse and D. R. Woodall, Choosability of powers of circuits, Graphs and Combinatorics 19 (2003), 137–-144.

\bibitem{RWest}
R. Ramamurthi and D. B. West, Hypergraph extension of the Alon-Tarsi list coloring theorem, Combinatorica 25(3) (2005), 355--366.

\bibitem{Ryser}
H.\ J.\ Ryser, Combinatorial Mathematics, The Carus mathematical monographs series, published by The Mathematical Association of America, 1963.

\bibitem{DES}
D. E. Scheim, The number of edge 3-colorings of a planar cubic
graph as a permanent, Discrete Math. 8 (1974), 377--382.

\bibitem{Tarsi}
M. Tarsi, The graph polynomial and the number of proper vertex colorings, Ann. Inst. Fourier, Grenoble \textbf{49}, 3 (1999), 1089--1093.

\bibitem{Tes1}
B. A. Tesman, List $T$-colorings of graphs, Discrete Applied Mathematics 45
(1993), 277--289.

\bibitem{Tes2}
B. A. Tesman, $T$-colorings, list $T$-colorings, and set $T$-colorings of graphs, manuscript.

\bibitem{Tru}
M. Truszczy\'nski, Results on uniquely colorable graphs, Colloquia Math. Soc. J\'anos Bolyai 37 (1981), 733--746. 

\bibitem{Vizing}
V. G. Vizing, Coloring the vertices of a graph in prescribed colors (in Russian),
Diskret. Analiz. No. 29, Metody Diskret. Anal. v. Teorii Kodov i Shem 101 (1976),
3–-10.

\bibitem{West}
D.\ B.\ West, {\bf Introduction to Graph Theory}, Prentice Hall,
2001.

\bibitem{Xu}
S. Xu, The size of uniquely colorable graphs, J. Combinatorial Theory Ser B 50 (1990), 319--320.

\end{thebibliography}
\end{document}